\documentstyle[amscd,amssymb,verbatim,12pt]{amsart}
\input diagrams
\diagramstyle[PostScript=dvips]

\newcommand{\cohproj}{\operatorname{cohproj}}

\newcommand{\can}{\operatorname{can}}

\newcommand{\tr}{\operatorname{tr}}

\newcommand{\coker}{\operatorname{coker}}

\renewcommand{\Re}{\operatorname{Re}}
\newcommand{\SL}{\operatorname{SL}}

\newcommand{\lan}{\langle}
\newcommand{\ran}{\rangle}
\newcommand{\Coh}{\operatorname{Coh}}

\newcommand{\CC}{{\cal C}}

\newcommand{\Proj}{\operatorname{Proj}}

\newcommand{\si}{\sigma}

\newcommand{\de}{\delta}
\newcommand{\eps}{\epsilon}

\renewcommand{\ker}{\operatorname{ker}}

\numberwithin{equation}{section}

\newtheorem{thm}{Theorem}[section]
\newtheorem{prop}[thm]{Proposition}
\newtheorem{lem}[thm]{Lemma}
\newtheorem{cor}[thm]{Corollary}
\newenvironment{rem}{\vspace{3mm}\noindent
{\bf Remark.}}{\vspace{3mm}}
\newenvironment{defi}{\vspace{3mm}\noindent
{\bf Definition.}}{\vspace{3mm}}

\newcommand{\Pf}{\noindent {\it Proof}}
\newcommand{\id}{\operatorname{id}}

\newcommand{\ov}{\overline}

\renewcommand{\Im}{\operatorname{Im}}

\newcommand{\rk}{\operatorname{rk}}
\newcommand{\ra}{\rightarrow}
\renewcommand{\AA}{{\cal A}}
\newcommand{\FF}{{\cal F}}

\newcommand{\SS}{{\cal S}}

\newcommand{\Hom}{\operatorname{Hom}}
\newcommand{\Ext}{\operatorname{Ext}}
\newcommand{\End}{\operatorname{End}}

\renewcommand{\a}{\alpha}

\newcommand{\la}{\lambda}
\newcommand{\th}{\theta}
\newcommand{\C}{{\Bbb C}}
\newcommand{\R}{{\Bbb R}}
\newcommand{\Z}{{\Bbb Z}}
\newcommand{\Q}{{\Bbb Q}}

\newcommand{\wt}{\widetilde}

\newcommand{\sign}{\operatorname{sign}}
\newcommand{\sub}{\subset}
\newcommand{\ed}{\qed\vspace{3mm}}

\title{Classification of holomorphic vector bundles on noncommutative two-tori}
\author{A. Polishchuk}
\thanks{This work was partially supported by NSF grant DMS-0302215}

\begin{document}
\begin{abstract} We prove that every holomorphic vector bundle on a noncommutative two-torus $T$
can be obtained by successive extensions from standard holomorphic bundles considered in \cite{PS}.
This implies that the category of holomorphic bundles on $T$ is equivalent to the heart of
certain $t$-structure on the derived category of coherent sheaves on an elliptic curve.
\end{abstract}
\maketitle

\bigskip

\section{Introduction}

In this paper we continue the study of holomorphic bundles on noncommutative two-tori that was begun in \cite{PS}.
Recall that for every $\th\in\R\setminus\Q$ and $\tau\in\C\setminus \R$ we considered in \cite{PS} holomorphic vector bundles
on a noncommutative complex torus $T=T_{\th,\tau}$. By definition, the algebra $A_{\th}$ of smooth
functions on $T$ consists of series $\sum_{(m,n)\in\Z^2} a_{m,n} U_1^m U_2^n$ where the coefficients
$a_{m,n}\in\C$ rapidly decrease at infinity and the multiplication is defined using the rule
$$U_1U_2=\exp(2\pi i \th)U_2U_1.$$
We consider the derivation $\de=\de_{\tau}:A_{\th}\ra A_{\th}$ defined by
$$\de(\sum a_{m,n} U_1^m U_2^n)=2\pi i\sum_{m,n} (m\tau+n)a_{m,n}U_1^mU_2^n$$
as an analogue of the $\ov{\partial}$-operator. A holomorphic bundle over $T$ is a pair
$(E,\ov{\nabla})$ consisting of a finitely generated projective right $A_{\th}$-module $E$ and
an operator $\ov{\nabla}:E\ra E$ satisfying the Leibnitz identity
$$\ov{\nabla}(ea)=\ov{\nabla}(e)a+e\de(a),$$
where $e\in E$, $a\in A_{\th}$. There is an obvious definition of a holomorphic map
between holomorphic bundles, so we can define the category $\CC(T)$ of holomorphic bundles
on $T$.

For every pair of relatively prime integers $(c,d)$ such that $c\th+d>0$ and a complex number
$z$ we define a {\it standard holomorphic bundle} $(E_{d,c}(\th),\ov{\nabla}_z)$ as follows.
If $c\neq 0$ then 
$$E_{d,c}(\th)=\SS(\R\times\Z/c\Z)=\SS(\R)^{|c|},$$ 
where $\SS(\R)$ is the Schwartz space of functions on $\R$,
with the $A_{\th}$-action defined by 
$$fU_1(x,\a)=f(x-\frac{1}{\mu},\a-1),\ fU_2(x,\a)=\exp(2\pi i(x-\frac{\a d}{c}))f(x,\a),
$$
where $x\in\R$, $\a\in\Z/c\Z$, $\mu=\frac{c}{c\th+d}$. 
The operator $\ov{\nabla}_z$ on this space is given by
\begin{equation}\label{nablaeq}
\ov{\nabla}_z(f)=\frac{\partial f}{\partial x}+2\pi i(\tau\mu x+z)f.
\end{equation}
For $c=0$ and $d=1$ we set $E_{1,0}(\th)=A_{\th}$ with the natural right $A_{\th}$-action
and the operator $\ov{\nabla}_z$ is given by
$$\ov{\nabla}_z(a)=\de(a)+2\pi i z a.$$
We define degree, rank and slope of a bundle $E=E_{d,c}(\th)$ by setting
$\deg(E)=c$, $\rk(E)=c\th+d$ and $\mu(E)=\deg(E)/\rk(E)$.
Note that $\rk(E)>0$ and $\mu=\mu(E)$ in the formulae above.

According to the theorem of Rieffel (see \cite{R}) every finitely generated projective
right $A_{\th}$-module is isomorphic to $E=E_{d,c}(\th)^{\oplus n}$ for some $(c,d)$ as above
and $n\ge 0$. Moreover, the degree and rank defined above extend to additive functions on the
category of finitely generated projective $A_{\th}$-modules.

The category of holomorphic bundles $\CC=\CC(T)$ has a natural
structure of a $\C$-linear exact category. In particular, for every pair of holomorphic
bundles $E_1$ and $E_2$ we can form the vector space $\Ext^1_{\CC}(E_1,E_2)$ parametrizing
extensions of $E_1$ by $E_2$. Sometimes we will also use the notation $\Ext^0_{\CC}:=\Hom_{\CC}$.
Let $\CC'\sub\CC$ be the minimal full subcategory of $\CC$ containing all standard holomorphic bundles
and closed under extensions. Our main result is the following theorem.

\begin{thm}\label{mainthm} One has $\CC'=\CC$.
\end{thm}

Combining this theorem with the study of the category $\CC'$ in \cite{PS} we obtain the following result.

\begin{cor} The category $\CC$ is abelian. It is equivalent to the heart $\CC^{\th}$ of the $t$-structure on
the derived category of coherent sheaves on the elliptic curve $\C/\Z+\Z\tau$, associated with $\th$
(see section 3 of \cite{PS} or section \ref{ampleseq} below).
\end{cor}

\begin{rem} Recall that we always assume $\th$ to be irrational. For rational $\th$
the category $\CC$ will not be abelian.
\end{rem}

\begin{cor} For every indecomposable
holomorphic bundle $E$ on $T$ there exists a standard holomorphic bundle $\ov{E}$ and
a filtration $0=E_0\sub E_1\sub\ldots\sub E_n=E$ by holomorphic subbundles such
that all quotients $E_i/E_{i-1}$ are isomorphic to $\ov{E}$.
\end{cor}

The proof of Theorem \ref{mainthm} consists of two steps. First, we develop the cohomology theory
for holomorphic bundles on $T$ and prove the analogues of the standard theorems for them
(such as finiteness, Riemann-Roch and Serre duality). Then we combine
these results with the techniques of \cite{P-real}
where the category $\CC'$ was described in terms of coherent modules over certain algebra.

\noindent
{\it Acknowledgments}. 
Parts of this paper were written during the author's visits to Max-Planck-Institut
f\"ur Mathematik in Bonn and the Institut des Hautes \'Etudes Scientifiques.
I'd like to thank these institutions for hospitality and support. 

\section{Cohomology of holomorphic bundles on noncommutative two-tori}

\subsection{Cohomology and $\Ext$-spaces}

Let $(E,\ov{\nabla})$ be a holomorphic bundle on $T=T_{\th,\tau}$. Then the {\it cohomology} of $E$
is defined by
$$H^i(E)=H^i(E,\ov{\nabla})=H^i(E\stackrel{\ov{\nabla}}{\ra}E)$$
where $i=0$ or $i=1$. Thus, $H^0(E)=\ker(\ov{\nabla})$, $H^1(E)=\coker(\ov{\nabla})$.
These spaces are closely related to $\Ext^i$-spaces in the category of holomorphic
bundles (where $i=0$ or $i=1$). 
To explain this connection we have to use Morita equivalences between noncommutative tori.
Recall that for every standard bundle $E_0=E_{d,c}(\th)$ the algebra of endomorphisms
$\End_{A_{\th}}(E_0)$ can be identified with the algebra $A_{\th'}$ for some $\th'\in\R$.
In fact, $\th'=\frac{a\th+b}{c\th+d}$, where $a$ and $b$ are chosen in such a way that
$\left(\matrix a & b \\ c & d\endmatrix\right)\in\SL_2(\Z)$.
Furthermore, if $E_0$ is equipped with a standard holomorphic structere $\ov{\nabla}$ then
the formula $\phi\mapsto [\ov{\nabla},\phi]$ defines a derivation of $\End_{A_{\th}}(E_0)\simeq A_{\th'}$,
hence the corresponding torus $T_{\th'}$ is equipped with a complex structure. In fact,
this derivation on $A_{\th'}$ is equal to $\de_{\tau}/\rk(E_0)$, where $\tau$ is the same parameter
that was used to define the complex structure on $T_{\th}$ (see Proposition 2.1 of \cite{PS}). 
Now one can define the Morita equivalence 
$$\CC(T_{\th',\tau})\ra\CC(T_{\th,\tau}):E\mapsto E\otimes_{A_{\th'}}E_0,$$
where the tensor product is equipped with the complex structure
$$\ov{\nabla}(e\otimes e_0)=\frac{1}{\rk(E_0)}\ov{\nabla}_E(e)\otimes e_0+e\otimes\ov{\nabla}_{E_0}(e_0)$$ 
(see Propositions 2.1 and 3.2 of \cite{PS}). This functor sends standard holomorphic bundles on
$T_{\th',\tau}$ to standard holomorphic bundles on $T_{\th,\tau}$. The inverse functor is 
\begin{equation}\label{dualfuneq}
\CC(T_{\th,\tau})\ra\CC(T_{\th',\tau}):E\mapsto \Hom_{A_{\th}}(E_0,E),
\end{equation}
where the latter space has a natural right action of $A_{\th'}\simeq \End_{A_{\th}}(E_0)$.
Now we can formulate the connection between the cohomology and $\Ext$-groups. For every
holomorphic bundle $E$ and a standard holomorphic bundle $E_0$ on $T=T_{\th,\tau}$ one has
a natural isomorphism
\begin{equation}\label{homeq1}
\Ext^i_{\CC}(E_0,E)\simeq H^i(\Hom_{A_{\th}}(E_0,E)),
\end{equation}
where $\Hom_{A_{\th}}(E_0,E)$ is viewed as a holomorphic bundle on $T_{\th'}$ (the proof is
similar to Proposition 2.4 of \cite{PS}).
Note that for an arbitrary pair of holomorphic bundles $E_1$ and $E_2$ one can still
define an operator $\ov{\nabla}$ on $\Hom_{A_{\th}}(E_1,E_2)$ such that the
analogue of isomorphism (\ref{homeq1}) holds. However, we have a natural interpretation
of $\Hom_{A_{\th}}(E_1,E_2)$ as a holomorphic bundle on some noncommutative two-torus only in the case when
one of the bundles $E_1$ or $E_2$ is standard (see (\ref{homeq2}) below).

\subsection{Duality and metrics}\label{dualsec}

One can define the category of {\it left holomorphic bundles} on $T$ by replacing right $A_{\th}$-modules
with left ones and changing the Leibnitz identity appropriately. The definition of cohomology for these bundles remains the same.
There is a natural duality functor $E\mapsto E^{\vee}$ that associates
to a (right) holomorphic bundle $E$ the left holomorphic bundle $\Hom_{A_{\th}}(E,A_{\th})$.
More generally, for every standard holomorphic bundle $E_0$ we can consider
$\Hom_{A_{\th}}(E,E_0)$ as a left module over $\End_{A_{\th}}(E_0)\simeq A_{\th'}$ equipped with an
induced holomorphic structure. Then the natural isomorphism
\begin{equation}\label{homeq2}
\Ext^i_{\CC}(E,E_0)\simeq H^i(\Hom_{A_{\th}}(E,E_0))
\end{equation}
allows to view $\Ext^i_{\CC}(E,E_0)$ as cohomology of a holomorphic bundle on $T_{\th'}$.
Using duality the functor (\ref{dualfuneq}) for a standard holomorphic bundle $E_0$
can be rewritten as the usual Morita functor due to the isomorphism
\begin{equation}\label{tensorhom}
\Hom_{A_{\th}}(E_0,E)\simeq E\otimes_{A_{\th}} E_0^{\vee}.
\end{equation}

For a standard holomorphic bundle $E_0$ on $T_{\th,\tau}$ the dual bundle $E_0^{\vee}$ can also be considered as a right 
holomorphic bundle on $T_{\th',\tau}$, where $\End_{A_{\th}}(E_0)\simeq A_{\th'}$. In fact, it is again a standard holomorphic 
bundle (see Corollary 2.3 of \cite{PS}). More precisely, for $E_0=E_{d,c}(\th)$ we have $\th'=\frac{a\th+b}{c\th+d}$, where
$\left(\matrix a & b \\ c & d\endmatrix\right)\in\SL_2(\Z)$. 
The left action of $A_{\th'}$ on $E_{d,c}(\th)=\SS(\R\times\Z/c\Z)$ (where $c\neq 0$) is defined by the formulae
$$U_1f(x,\a)=f(x-\frac{1}{c},\a-a),\ 
U_2f(x,\a)=\exp(2\pi i(\frac{x}{c\th+d}-\frac{\a}{c}))f(x,\a),$$
where $x\in\R$, $\a\in\Z/c\Z$.
We can identify $E_{d,c}(\th)^{\vee}$ considered as a right $A_{\th'}$-module
with $E_{a,-c}(\th')$ using the natural pairing
$$t:E_{a,-c}(\th')\otimes E_{d,c}(\th)\ra A_{\th}$$
constructed as follows (see Proposition 1.2 of \cite{PS}). First, we define the map
$$b:E_{a,-c}(\th')\otimes E_{d,c}(\th)\ra\C$$
by the formula
$$b(f_1,f_2)=\sum_{\a\in\Z/c\Z}\int_{x\in\R} f_1(\frac{x}{c\th+d},\a)f_2(x,-a\a)dx.$$
Then $t$ is given by
$$t(f_1,f_2)=\sum_{(m,n)\in\Z^2}
U_1^{m}U_2^{n}b(U_2^{-n}U_1^{-m}f_1\otimes f_2).$$ 
The corresponding isomorphism 
$$E_{d,c}(\th)^{\vee}\simeq E_{a,-c}(\th')$$
is compatible with the $A_{\th}-A_{\th'}$-bimodule structures and with holomorphic structures
(see Corollary 2.3 of \cite{PS}).
Note that $b=\tr\circ t$, where $\tr:A_{\th}\ra\C$ is the trace
functional sending $\sum a_{m,n} U_1^m U_2^n$ to $a_{0,0}$.

On the other hand, we can define a $\C$-antilinear isomorphism 
$$\si:E_{d,c}(\th)\ra E_{a,-c}(\th')$$
by the formula
$$\si(f)(x,\a)=\ov{f((c\th+d)x,-a\a)}.$$ 
This isomorphism satisfies 
\begin{equation}\label{sigmaeq}
\si(bea)=a^*\si(e)b^*
\end{equation}
for $e\in E_{d,c}(\th)$,
$a\in A_{\th}$, $b\in A_{\th'}$, where $*:A_{\th}\ra A_{\th}$ is the $\C$-antilinear anti-involution sending
$U_i$ to $U_i^{-1}$.
In view of the identification of $E_{a,-c}(\th')$ with the dual bundle to $E_{d,c}(\th)$ the isomorphism $\si$
should be considered as an analogue of the Hermitian metric on $E_{d,c}(\th)$. 
The corresponding analogue of the scalar product on global sections is simply the Hermitian form on $E_{d,c}(\th)$
given by the formula
\begin{equation}\label{hermform}
\lan f_1,f_2\ran=b(\si(f_2),f_1)=\sum_{\a\in\Z/c\Z}\int_{x\in\R}f_1(x,\a)\ov{f_2(x,\a)}dx,
\end{equation}
where $f_1,f_2\in E_{d,c}(\th)$.
We can also define the corresponding $L^2$-norm: $||f||_0^2=\lan f,f\ran$.
The above Hermitian form is related to the structure of $A_{\th'}-A_{\th}$-bimodule on $E_{d,c}(\th)$ in the following way:
$$\lan f_1, af_2b\ran=\lan a^*f_1b^*,f_2\ran,$$
where $a\in A_{\th'}$, $b\in A_{\th}$ (this is a consequence of (\ref{sigmaeq}) and of Lemma 1.1 of \cite{PS}).

In the case of the trivial bundle $E_{1,0}(\th)=A_{\th}$ we can easily modify the above definitions. First of all,
$\th'=\th$ and the dual bundle is still $A_{\th}$. The role of $\si$ is played by $*:A_{\th}\ra A_{\th}$ and
the Hermitian form on $A_{\th}$ is given by $\lan a,b\ran=\tr(ab^*)$. The corresponding $L^2$-norm is   
$$||\sum a_{m,n} U_1^m U_2^n||^2_0=\sum |a_{m,n}|^2.$$

Note that the operator $\ov{\nabla}_z$ on $E_{d,c}(\th)$ admits an adjoint operator $\ov{\nabla}_z^*$ with
respect to the Hermitian metrics introduced above. Namely, for $c\neq 0$ it is given by
$$\ov{\nabla}_z^*(f)=-\frac{\partial f}{\partial x}-2\pi i(\ov{\tau}\mu x+\ov{z})f,$$
while for $c=0$ we have $\ov{\nabla}_z^*=-\de_{\ov{\tau}}-2\pi i z\id$ on $A_{\th}$.
In either case we have 
\begin{equation}\label{curveq}
\ov{\nabla}_z\ov{\nabla}_z^*-\ov{\nabla}_z^*\ov{\nabla}_z=\la\cdot\id
\end{equation}
for some constant $\la\in\R$.

It follows that for an arbitrary holomorphic structure $\ov{\nabla}$ on
$E=E_{d,c}(\th)^{\oplus n}$ there exists an adjoint operator $\ov{\nabla}^*:E\ra E$ with respect to the
above Hermitian metric. Indeed, we can write $\ov{\nabla}=\ov{\nabla}_0+\phi$, where $\ov{\nabla}_0$
is the standard holomorphic structure and set $\ov{\nabla}^*=\ov{\nabla}_0^*+\phi^*$.

\subsection{Sobolev spaces}

The idea to consider Sobolev spaces for bundles on noncommutative tori is due to M.~Spera (see \cite{Spera}, \cite{Spera2}).
Let $(E,\ov{\nabla})$ be a standard holomorphic bundle on $T_{\th,\tau}$.
For $s\in\Z$, $s\ge 0$ we define the $s$-th Sobolev norm on $E$ by setting
$$||e||_s^2=\sum_{i=0}^s ||\ov{\nabla}^i e||_0^2,$$
where $||e||_0$ is the $L^2$-norm on $E$.
We define $W_s(E)$ to be the completion of $E$ with respect to this norm.
Note that there is a natural embedding $W_{s+1}(E)\sub W_s(E)$.
We can define analogous spaces for $E^{\oplus n}$ in an obvious way.

All the definitions above make sense also for rational $\th$.
Moreover, for $\th\in\Z$ the space $E$ can be identified with the space of smooth section of a holomorphic
vector bundle $V$ on an elliptic curve $\C/\Z+\Z\tau$ in such a way that
$\ov{\nabla}$ corresponds to the $\ov{\partial}$-operator. Furthermore, the $L^2$-norm above
corresponds to the $L^2$-norm with respect to a Hermitian metric on $V$ that has constant curvature.
This implies that in this case the Sobolev spaces $W_s(E)$ coincide with the corresponding Sobolev spaces constructed for
the holomorphic bundle $V$. Indeed, using the equation (\ref{curveq}) it is easy to see
that the norm $||e||_s$ is equivalent to the norm given by 
$$(||e||'_s)^2=\sum_{i=0}^s\lan e, (\ov{\nabla}^*\ov{\nabla})^i e\ran$$
which is equivalent to the standard Sobolev norm.

An important observation is that the operator $\ov{\nabla}_z$ defined by (\ref{nablaeq})
depends only on $\tau\mu$, where $\mu$ is the slope of the bundle, 
so it is the same for the bundle $E_{d,c}(\th)$ on $T_{\th,\tau}$
and the bundle $E_{d,c}(\sign(c)N)$ on the commutative torus $T_{\sign(c)N,\tau'}$,
where $N$ is a large enough integer so that $|c|N+d>0$, $\tau'=(|c|N+d)/(c\th+d)$.
Therefore, the sequences of spaces $(W_s(E))$ in these two cases are the same.
Hence, the following standard results about Sobolev spaces in the commutative case
extend immediately to our situation (the first two are analogues of Rellich's lemma and Sobolev's lemma).
In all these results $E$ is a direct sum of a finite number of copies of a standard holomorphic
bundle.

\begin{lem}\label{Rellichlem} 
The embedding $W_s(E)\sub W_{s-1}(E)$ is a compact operator.
\end{lem}

\begin{lem}\label{Sobolevlem} 
One has $E=\cap_{s\ge 0} W_s(E)$.
\end{lem}

\begin{lem}\label{Sobdiflem} 
The operator $\ov{\nabla}$ extends to a bounded operator
$W_s(E)\ra W_{s-1}(E)$
\end{lem}

The following result is the only noncommutative contribution to the techniques of Sobolev spaces, however,
it is quite easy.

\begin{lem}\label{Sobendlem} 
For every $\phi\in\End_{A_{\th}}(E)$ the operator $\phi:E\ra E$
extends to a bounded operator 
$W_s(E)\ra W_s(E)$.
\end{lem}

\Pf . It suffices to prove that for every $s\ge 0$ one has
$$||\phi e||_s\le C\cdot ||e||_s$$
for some constant $C>0$. By our assumption $E\simeq E_0^{\oplus N}$ for some standard bundle $E_0$. 
Identifying $\End_{A_{\th}}(E_0)$ with $A_{\th'}$ for some $\th'\in\R$ 
we can write $\phi=\sum a_{m,n} U_1^m U_2^n$ where 
$U_1$ and $U_2$ are unitary generators of $A_{\th'}$, $a_{m,n}$ are complex $N\times N$ matrices. Since
$U_1$ and $U_2$ act on $E$ by unitary operators, it follows that
$$||\phi e||_0\le C(\phi)\cdot||e||_0$$
for $e\in E$, where $C(\phi)=\sum_{m,n} ||a_{m,n}||$ (here $||a||$ denote the norm of a matrix $a$). 
Applying the Leibnitz rule repeatedly we derive similarly that
$$\sum_{i=0}^s||\ov{\nabla}^i(\phi e)||_0^2\le\sum_{i=0}^s c_i\cdot ||\ov{\nabla}^i e||_0^2$$
for some constants $c_i>0$ which implies the result.
\ed

It is convenient to extend the definition of the chain $\ldots \sub W_1(E)\sub W_0(S)$
to the chain of embedded spaces 
$$\ldots \sub W_1(E)\sub W_0(S)\sub W_{-1}(E)\sub\ldots$$
by setting $W_{-s}(E)=\ov{W_s(E)}^*$ (the space of $\C$-antilinear functionals) and using the natural
Hermitian form of $W_0(E)$. 
It is easy to see that the results of this section hold for all integer values of $s$.

\begin{lem}\label{adjointlem} Let $\ov{\nabla}:E\ra E$ be a (not necessarily standard) holomorphic structure on $E$.
Then the operators $\ov{\nabla}$ and $\ov{\nabla}^*$ can be extended to bounded operators $W_s(E)\ra W_{s-1}(E)$ for every $s\in\Z$.
\end{lem}

\Pf . Let $\ov{\nabla}_0$ be a standard holomorphic structure on $E$.
Then $\ov{\nabla}=\ov{\nabla}_0+\phi$ for some $\phi\in\End_{A_{\th}}(E)$. By Lemma \ref{Sobdiflem} (resp., Lemma
\ref{Sobendlem}) there exist a continuous extension $\ov{\nabla}_0:W_s(E)\ra W_{s-1}(E)$ (resp.,
$\phi:W_s(E)\ra W_s(E)$). Hence, $\ov{\nabla}$ extends to a family of continuous operators $\ov{\nabla}(s):W_s(E)\ra W_{s-1}(E)$ for $s\in\Z$.
The extensions of $\ov{\nabla}^*$ are given by the adjoint operators $\ov{\nabla}(-s+1)^*:W_s(E)\ra W_{s-1}(E)$.
\ed

\subsection{Applications to cohomology}

We begin our study of cohomology with standard holomorphic bundles.

\begin{prop}\label{Qprop} 
Let $(E,\ov{\nabla})$ be a direct sum of several copies of a standard holomorphic bundle on $T_{\th,\tau}$.

\noindent
(i) The cohomology spaces $H^0(E)$ and $H^1(E)$ are finite-dimensional and for $\Im(\tau)<0$ one has
$$\chi(E)=\dim H^0(E)-\dim H^1(E)=\deg(E).$$

\noindent
(ii) There exists an operator $Q:E\ra E$ such that 
$$\id-Q\ov{\nabla}=\pi_{\ker{\ov{\nabla}}},$$
$$\id-\ov{\nabla}Q=\pi_{\ov{\nabla}(E)^{\perp}},$$
where $\ov{\nabla}(E)^{\perp}\sub E$ is the orthogonal complement to $\ov{\nabla}(E)\sub E$,
for a finite-dimensional subspace $V\sub E$ we denote by $\pi_V:E\ra V$ the orthogonal projection.

\noindent
(iii) If $\Im(\tau)<0$ and $\deg(E)>0$ then $H^1(E)=0$.

\noindent
(iv) The operator $Q:E\ra E$ extends to a bounded operator  
$W_s(E)\ra W_{s+1}(E)$. 

\noindent
(v) For every $e\in W_0(E)$ one has
$$||Qe||_0\le \frac{1}{2\sqrt{\pi|\Im(\tau)\mu(E)|}}||e||_0.$$
\end{prop}

\Pf . In the commutative case the assertions (i)-(iii) are well known. For example, the operator $Q$ is given by 
of $\ov{\partial}^*G$, where $G$ is the Green operator for the $\ov{\partial}$-Laplacian. The condition
$\Im(\tau)<0$ corresponds to the way we define the operator $\de_{\tau}$ on $A_{\th}$ (see Proposition 3.1 of \cite{PS}).
As before we can deduce (i)-(iii) in general from the commutative case. One can also prove these assertions directly
in the noncommutative case (see Proposition 2.5 of \cite{PS} for the proofs of (i) and (iii)).
The assertion (iv) follows immediately from the identity 
$$\ov{\nabla}^nQ=\ov{\nabla}^{n-1}(\id-\pi_{\ov{\nabla}(E)^{\perp}}).$$
To prove (v) we can assume that $E=E_{d,c}(\th)$,
where $c\neq 0$ and $\ov{\nabla}=\ov{\nabla}_z$ for some $z\in\C$.
Then the space $W_0(E)$ is the orthogonal sum of $|c|$ copies of $L^2(\R)$. Moreover,
the operator $\ov{\nabla}$ respects this decomposition and restricts to the operator
$$f\mapsto f'+(ax+z)f$$
on each copy, where $a=2\pi i\tau \mu(E)$. Hence, the operator $Q$ also respects this decomposition and it 
suffices to consider its restriction to one copy of $L^2(\R)$.  
Since $\Re(a)\neq 0$, by making the unitary transformation
of the form $\wt{f}(x)=\exp(i tx)f(x+t')$ for some $t,t'\in\R$ we can reduce ourselves to the case $z=0$.
Furthermore, the transformation of the form
$\wt{f}=\exp(i \Im(a)x^2/2)f$ gives a unitary equivalence with the operator
$\ov{\nabla}:f\mapsto f'+\la xf$ where $\la=\Re(a)$.
Consider the following complete orthogonal system of functions in $L^2(\R)$: 
$$f_n(x)=H_n(\sqrt{|\la|}x)\exp(-|\la|\frac{x^2}{2}), \ n=0,1,2,\ldots,$$
where $H_n(x)=(-1)^n \exp(x^2)\frac{d^n}{dx^n}(\exp(-x^2))$ are Hermite polynomials
($(f_n)$ is an eigenbasis of the operator $f\mapsto -f''+\la^2x^2$). Note that
$$||f_n||_0^2=\frac{1}{\sqrt{\la}}\int_{\R}H_n(x)^2\exp(-x^2)dx=\frac{2^n\cdot n!\cdot\sqrt{\pi}}{\sqrt{\la}}.$$

Assume first that $\la>0$. Then using the formula $H_n'(x)=2n H_{n-1}(x)$ we obtain
$$\ov{\nabla}(f_n)=2n\sqrt{\la}f_{n-1}$$
for $n>0$ and $\ov{\nabla}(f_0)=0$.
Therefore, in this case
$$Q(f_n)=\frac{1}{(2n+2)\sqrt{\la}}f_{n+1}$$
for all $n\ge 0$. Hence,
$$\frac{||Q(f_n)||_0}{||f_n||_0}=\frac{||f_{n+1}||_0}{(2n+2)\sqrt{\la}||f_n||_0}=\frac{1}{\sqrt{(2n+2)\la}}\le\frac{1}{\sqrt{2\la}}$$
which implies (v) in this case. 

Now assume that $\la<0$. Then using the formula $H_n'(x)=2xH_n(x)-H_{n+1}(x)$ we find
$$\ov{\nabla}(f_n)=-\sqrt{|\la|}f_{n+1}.$$
Hence,
$$Q(f_n)=-\frac{1}{\sqrt{|\la|}f_{n-1}}$$
for $n>0$ and $Q(f_0)=0$. It follows
that 
$$\frac{||Q(f_n)||_0}{||f_n||_0}=\frac{||f_{n-1}||_0}{\sqrt{|\la|}||f_n||_0}=\frac{1}{\sqrt{2n|\la|}}\le
\frac{1}{\sqrt{2|\la|}}$$
for $n\ge 1$, which again implies our statement.
\ed

Now we are ready to prove results about cohomology of arbitrary holomorphic bundles. We will use
the following well known lemma.

\begin{lem}\label{compactlem} 
Let $L:W\ra W'$ and $L':W'\ra W$ be bounded operators between Banach spaces
such that $L'L=\id+C$, $LL'=\id+C'$ for some compact operators $C:W\ra W$ and $C':W'\ra W'$.
Then the operator $L$ is Fredholm.
\end{lem}

\begin{thm}\label{RRthm} 
(i) For every holomorphic bundle $(E,\ov{\nabla})$ on $T_{\th,\tau}$ the spaces
$H^0(E)$ and $H^1(E)$ are finite-dimensional. 

\noindent
(ii) If $\Im\tau<0$ then
$$\chi(E)=\dim H^0(E)-\dim H^1(E)=\deg(E).$$

\noindent
(iii) Let us equip $E$ with a metric by identifying it with the direct sum of several copies
of a standard bundle. Then one has the following orthogonal decompositions
$$E=\ker(\ov{\nabla})\oplus\ov{\nabla}^*(E),$$
$$E=\ker(\ov{\nabla}^*)\oplus\ov{\nabla}(E),$$
where $\ov{\nabla}^*:E\ra E$ is the adjoint operator to $\ov{\nabla}$.
\end{thm}

\Pf . Let us write the holomorphic structure on $E$ in the form
$$\ov{\nabla}=\ov{\nabla}_0+\phi,$$
where $(E,\ov{\nabla}_0)$ is holomorphically isomorphic to the direct sum of several copies of a standard holomorphic bundle,
$\phi\in\End_{A_{\th}}(E)$.
By Lemma \ref{adjointlem},
$\ov{\nabla}$ has a bounded extension to an operator $\ov{\nabla}:W_s(E)\ra W_{s-1}(E)$ for every $s\in\Z$.

Consider the operator $Q:E\ra E$ constructed in Proposition \ref{Qprop} for the holomorphic
structure $\ov{\nabla}_0$. Then $Q$ extends to a bounded operator
$W_s(E)\ra W_{s+1}(E)$ for every $s\in\Z$. We have 
$$Q\ov{\nabla}=Q\ov{\nabla}_0+Q\phi=\id-\pi_0+Q\phi,$$
where $\pi_0$ is the orthogonal projection to the finite-dimensional space $\ker(\ov{\nabla}_0)\sub E$.
Clearly, $\pi_0$ defines a bounded operator $W_0(E)\ra \ker(\ov{\nabla}_0)$.
Hence, the operator $C=Q\ov{\nabla}-\id:W_s(E)\ra W_s(E)$ factors as a composition
of some bounded operator $W_s(E)\ra W_{s+1}(E)$ with the embedding $W_{s+1}(E)\ra W_s(E)$.
By Lemma \ref{Rellichlem} this implies that $C$ is a compact operator. Similarly,
the operator $C'=\ov{\nabla}Q-\id:W_s(E)\ra W_s(E)$ is compact. Applying Lemma \ref{compactlem}
we deduce that $\ov{\nabla}:W_s(E)\ra W_{s-1}(E)$ is a Fredholm operator.
This immediately implies that $H^0(E)$ is finite-dimensional. Moreover, we claim that 
$$\ker(\ov{\nabla}:W_s(E)\ra W_{s-1}(E))=H^0(E)\sub E$$
for any $s\in\Z$.
Indeed, it suffices to check that if $\ov{\nabla}(e)=0$ for $e\in W_s(E)$ then
$e\in E$. Let us prove by induction in $t\ge s$ that $e\in W_t(E)$. Assume that
this is true for some $t$. Then 
$$e=Q\ov{\nabla}(e)-C(e)=-C(e)\in W_{t+1}(E).$$
Since $\cap_t W_t(E)=E$ by Lemma \ref{Sobolevlem} we conclude that $e\in E$.

Let $Q^*:W_s(E)\ra W_{s+1}(E)$ be the adjoint operator to $Q:W_{-s-1}(E)\ra W_{-s}(E)$, where $s\in\Z$.  
Then the operators $Q^*\ov{\nabla}^*-\id=(C')^*$ and $\ov{\nabla}^*Q^*-\id=C^*$ are compact. Thus, 
the same argument as before shows that for every $s\in\Z$ 
the operator $\ov{\nabla}^*:W_s(E)\ra W_{s-1}(E)$ is Fredholm and one has
$$\ker(\ov{\nabla}^*:W_s(E)\ra W_{s-1}(E))\sub E.$$ 
Next we claim that
$$E=\ker(\ov{\nabla}^*)\oplus \ov{\nabla}(E).$$
Since the orthogonal complement to $\ker(\ov{\nabla}^*)$ in $W_0(E)$ coincides
with the image of $\ov{\nabla}:W_1(E)\ra W_0(E)$, it suffices to prove
that $E\cap\ov{\nabla}(W_1(E))\sub \ov{\nabla}(E)$. But if $e=\ov{\nabla}(e_1)$
for some $e\in E$, $e_1\in W_1(E)$, then we can easily prove by induction in $s\ge 1$
that $e_1\in W_s(E)$. Indeed, assuming that $e_1\in W_s(E)$ we have
$$e_1=Q\ov{\nabla}(e_1)-C(e_1)=Q(e)-C(e_1)\in W_{s+1}(E).$$
A similar argument using the operator $Q^*$ shows that
$$E=\ker(\ov{\nabla})\oplus \ov{\nabla}^*(E).$$

Thus, we checked that $H^0(E)$ and $H^1(E)$ are finite-dimensional and that
$\chi(E)$ coincides with the index of the Fredholm operator 
$\ov{\nabla}=\ov{\nabla}_0+\phi:W_1(E)\ra W_0(E)$. Note that 
$$\ov{\nabla}_t:=\ov{\nabla}_0+t\phi:W_1(E)\ra W_0(E)$$
is a continuous family of Fredholm operators depending on $t\in [0,1]$. It follows that
the index of $\ov{\nabla}=\ov{\nabla}_1$ is equal to the index of $\ov{\nabla}_0$ computed in
Proposition \ref{Qprop}.
\ed

\begin{cor}\label{RRcor} 
If one of the holomorphic bundles $E_1$ and $E_2$ is standard then
the spaces $\Hom_{\CC}(E_1,E_2)$ and $\Ext^1_{\CC}(E_1,E_2)$ are finite-dimensional and
$$\chi(E_1,E_2):=\dim\Hom_{\CC}(E_1,E_2)-\dim\Ext^1_{\CC}(E_1,E_2)=
\rk(E_1)\deg(E_2)-\rk(E_2)\deg(E_1).$$
\end{cor}

The following vanishing result will play a crucial role in the proof of Theorem \ref{mainthm}.

\begin{thm}\label{vanishthm} 
Assume that $\Im(\tau)<0$.
For every holomorphic bundle $E$ on $T=T_{\th,\tau}$ there exists a constant $C=C(E)\in\R$ such that
for every standard holomorphic bundle $E_0$ on $T$ with $\mu(E_0)<C$
one has $\Ext^1_{\CC}(E_0,E)=0$. 
\end{thm}

\Pf . Let us choose a (non-holomorphic) isomorphism $E\simeq E_1^{\oplus N}$, where $E_1$ is a standard holomorphic
bundle on $T_{\th,\tau}$. Then we can write the holomorphic structure on $E$ as
$$\ov{\nabla}_E=\ov{\nabla}_0+\phi,$$
where $\ov{\nabla}_0$ comes from the standard holomorphic structure on $E_1$, $\phi\in\End_{A_{\th}}(E)$. 
Then for every standard holomorphic bundle $E_0$ we can consider the holomorphic bundle
$E'=\Hom_{A_{\th}}(E_0,E)$ on $T_{\th',\tau}$, where $\End_{A_{\th}}(E_0)=A_{\th'}$.
Note that $\Ext^1_{\CC}(E_0,E)\simeq H^1(E')$, so we want to prove that the latter group vanishes for $\mu(E_0)<<0$.
Recall that the holomorphic structure $\ov{\nabla}'$ on $E'$ is given by
$$\ov{\nabla}'(f)(e_0)=\rk(E_0)\cdot[\ov{\nabla}(f(e_0))-f(\ov{\nabla}_{E_0}(e_0))],$$
where $e_0\in E_0$, $f\in E'$ (see section 2.2 of \cite{PS}). 
The isomorphism $E\simeq E_1^{\oplus N}$ induces an isomorphism
$E'\simeq (E'_1)^{\oplus N}$, where $E'_1$ is the standard bundle $\Hom_{A_{\th}}(E_0,E_1)$ on $T_{\th',\tau}$.
Therefore, we have
$$\ov{\nabla}'=\ov{\nabla}'_0+\rk(E_0)\phi,$$ 
where $\ov{\nabla}'_0$ corresponds to the standard
holomorphic structure on $(E'_1)^{\oplus N}$ and $\phi$ is now considered as an $A_{\th'}$-linear endomorphism of $E'$. 
Note that by Proposition \ref{Qprop}(iii) we have 
$H^1(E',\ov{\nabla}'_0)=0$ as long as $\mu(E')>0$.
It is easy to compute that
$$\rk(E')=\rk(E)/\rk(E_0), \ \deg(E')=\rk(E)\rk(E_0)(\mu(E)-\mu(E_0)),$$
hence
$$\mu(E')=\rk(E_0)^2(\mu(E)-\mu(E_0)).$$
Therefore, $\mu(E')>0$ provided that $\mu(E_0)<\mu(E)$. In this case the operator $Q$ on $E'$ constructed in
Proposition \ref{Qprop} for the standard holomorphic structure $\ov{\nabla}'_0$ satisfies $\ov{\nabla}'_0Q=\id$.
Let $C_0$ be a constant such that
$$||Qe'||_0\le \frac{C_0}{\sqrt{\mu(E')}}||e'||_0$$
for $e'\in W_0(E')$ (see Proposition \ref{Qprop}(v)). 
Also, let us write $\phi\in\End_{A_{\th}}(E)\simeq\End_{A_{\th}}(E_1^{\oplus N})$ in the form
$\phi=\sum a_{m,n} U_1^m U_2^n$, where $U_1$ and $U_2$ are unitary generators of $\End_{A_{\th}}(E_1)$
and $a_{m,n}$ are $N\times N$ complex matrices.
Then we set $C(\phi)=\sum_{m,n} ||a_{m,n}||$ (the sum of the matrix norms of all coefficients).
Now we choose the constant $C<\mu(E)$ in such a way that
$$\frac{C_0C(\phi)}{\sqrt{\mu(E)-C}}<1.$$
Then for $\mu(E_0)<C$ we will have
$$\rk(E_0)\cdot ||\phi Qe'||_0\le \frac{\rk(E_0)C_0C(\phi)}{\sqrt{\mu(E')}}||e'||_0<\frac{C_0C(\phi)}{\sqrt{\mu(E)-C}}||e'||_0<r\cdot ||e'||_0$$
for some $0<r<1$. It follows from the above estimate that the operator $\id+\rk(E_0)\phi Q:W_0(E')\ra W_0(E')$ is invertible.
Therefore, we can define the operator
$$\wt{Q}=Q(\id+\rk(E_0)\phi Q)^{-1}:W_0(E')\ra W_1(E')$$
that satisfies 
$$(\ov{\nabla}'_0+\rk(E_0)\phi)\wt{Q}=\id.$$ 
Hence, the operator $\ov{\nabla}'=\ov{\nabla}'_0+\rk(E_0)\phi:W_1(E')\ra W_0(E')$ is surjective. But $H^1(E')$ can
be identified with the cokernel of this operator (see the proof of Theorem \ref{RRthm}),
so $H^1(E')=0$.
\ed

\subsection{Serre duality}

For every holomorphic bundle $E$ we have a natural pairing
$$E\otimes_{A_{\th}} E^{\vee}\ra A_{\th}:e\otimes f\mapsto f(e).$$
It is compatible with the $\ov{\partial}$-operators, so
it induces a pairing
\begin{equation}\label{Serpar}
H^{1-i}(E)\otimes H^i(E^{\vee})\ra H^1(A_{\th})\simeq\C
\end{equation}
for $i=0,1$.

\begin{thm}\label{Serthm} 
The pairing (\ref{Serpar}) is perfect.
\end{thm}

\Pf . Since we can switch $E$ and $E^{\vee}$, it suffices to consider the case $i=0$.
We choose an isomorphism of $E$ with a direct sum of several copies of a standard holomorphic
bundle $E_0$, so that we can talk about standard metrics and Sobolev spaces.
Note that the isomorphism $H^1(A_{\th})$ is induced by the trace functional $\tr:A_{\th}\ra\C:\sum a_{m,n}U_1^mU_2^n\mapsto a_{0,0}$.
Hence, the pairing (\ref{Serpar}) is induced by the pairing
$$b:E\otimes E^{\vee}\ra\C:e\otimes f\mapsto \tr(f(e))$$
that satisfies the identity  
\begin{equation}\label{beq}
\rk(E_0)b(\ov{\nabla}_E(e),e^{\vee})+b(e,\ov{\nabla}_{E^{\vee}}(e^{\vee}))=0,
\end{equation}
where $e\in E$, $e^{\vee}\in E^{\vee}$ (see Proposition 2.2 of \cite{PS}).
According to Theorem \ref{RRthm}(iii) we have orthogonal decompositions
$$E^{\vee}=\ker(\ov{\nabla}_{E^{\vee}})\oplus \ov{\nabla}^*_{E^{\vee}}(E^{\vee}),$$
$$E=\ker(\ov{\nabla}^*_E)\oplus\ov{\nabla}_E(E).$$
Therefore, it suffices to check that $b$ induces a perfect pairing between
$\ker(\ov{\nabla}^*_E)$ and $\ker(\ov{\nabla}_{E^{\vee}})$.
Let $\si:E\ra E^{\vee}$ be the $\C$-antilinear isomorphism defined in section \ref{dualsec}.
We claim that $\si$ maps $\ker(\ov{\nabla}^*_E)$ isomorphically onto $\ker(\ov{\nabla}_{E^{\vee}})$.
Since $b(e_1,\si(e_2))=\lan e_1,e_2\ran$ for $e_1,e_2\in E$, the theorem would immediately follow this.
To prove the claim it is enough to check that $\ov{\nabla}_{E^{\vee}}=-\rk(E_0)\si\ov{\nabla}_E^*\si^{-1}$.
To this end let us rewrite (\ref{beq}) as follows:
$$\rk(E_0)\lan \ov{\nabla}_E(e), \si^{-1} e^{\vee}\ran=-\lan e, \si^{-1}\ov{\nabla}_{E^{\vee}}(e^{\vee})\ran.$$
Since the left-hand side is equal to $\rk(E_0)\lan e,\ov{\nabla}_E^*\si^{-1}e^{\vee}\ran$ we conclude
that $\rk(E_0)\ov{\nabla}_E^*\si^{-1}=-\si^{-1}\ov{\nabla}_{E^{\vee}}$ as required.
\ed

Recall that we denote by $\CC'\sub\CC$ the full subcategory consisting of all successive extensions of 
standard holomorphic bundles. Theorem 3.8 of \cite{PS} implies that the derived category of $\CC'$ is
equivalent to the derived category of coherent sheaves on an elliptic curve. Therefore, the standard
Serre duality gives a functorial isomorphism
$$\Ext^1_{\CC'}(E_1,E_2)\simeq\Hom_{\CC'}(E_2,E_1)^*$$
for $E_1,E_2\in\CC'$. Note that we can replace here $\Ext^i_{\CC'}$ with $\Ext^i_{\CC}$.
Now using the above theorem we can extend this isomorphism to the case when only one of the objects $E_1,E_2$
belongs to $\CC'$.

\begin{cor}\label{Sercor} 
For every holomorphic bundles $E$ and $E_0$ such that $E_0\in\CC'$ the
natural pairings
$$\Ext^i_{\CC}(E_0,E)\otimes \Ext^{1-i}_{\CC}(E,E_0)\ra \Ext^1_{\CC}(E_0,E_0)\ra\C$$
for $i=0,1$ are perfect. Here the functional $\Ext^1_{\CC}(E_0,E_0)\ra\C$ is induced by
the Serre duality on $\CC'$. Therefore, we have functorial isomorphisms
\begin{equation}\label{Serisom}
\Ext^{1-i}_{\CC}(E,E_0)\wt{\ra} \Ext^i_{\CC}(E_0,E)^*
\end{equation}
for $E\in\CC$, $E_0\in\CC'$.
\end{cor}

\Pf . If $E_0$ is standard the assertion follows from Theorem \ref{Serthm}. It remains to observe that if for fixed $E\in\CC$
the map (\ref{Serisom}) is an isomorphism for some $E_0,E'_0\in\CC'$ then it is also an isomorphism for any extension of
$E_0$ by $E'_0$.
\ed

\section{Ampleness}
\label{amplesec}

\subsection{Ample sequences of standard holomorphic bundles}\label{ampleseq}

Let us start by recalling some basic notions concerning ample sequences in abelian categories and
associated $\Z$-algebras. The reader can consult \cite{P-coh} for more details.

\begin{defi}(see \cite{SVdB},\cite{P-coh}):
A sequence of objects $(E_n)_{n\in\Z}$ in a $\C$-linear abelian category $\AA$ is called {\it ample}
if the following two conditions hold:

\noindent
(i) for every surjection $E\ra E'$ in $\AA$ the induced map $\Hom(E_n,E)\ra \Hom(E_n,E')$ is surjective for
all $n<<0$; 

\noindent
(ii) for every object $E\in\AA$ and every $N\in\Z$ there exists a surjection $\oplus_{i=1}^s E_{n_i}\ra E$
where $n_i<N$ for all $i$.
\end{defi}

To a sequence $(E_n)_{n\in\Z}$ one can associate a so called $\Z$-algebra $A=\oplus_{i\le j}A_{ij}$, where $A_{ii}=\C$,
$A_{ij}=\Hom(E_i,E_j)$ for $i<j$, the multiplications $A_{jk}\otimes A_{ij}\ra A_{ik}$ are
induced by the composition in $\CC$.
One can define for $\Z$-algebras all the standard notions associated with graded algebras (see \cite{P-coh}).
In particular, we can talk about right $A$-modules: these have form $M=\oplus_{i\in\Z} M_i$ and the right $A$-action
is given by the maps $M_j\otimes A_{ij}\ra A_i$. The analogues of free $A$-modules are direct sums of the modules 
$P_n$, $n\in\Z$, defined by $(P_n)_i=A_{ni}$. We say that an $A$-module $M$ is {\it finitely generated} if
there exists a surjection $\oplus_{i=1}^s P_{n_i}\ra M$. A finitely generated $A$-module $M$ is called
{\it coherent} if for every morphism $f:P\ra M$, where $P$ is a finitely generated free module, the module
$\ker(f)$ is finitely generated. Finally, a $\Z$-algebra $A$ is called {\it coherent} if all the modules
$P_n$ are coherent and in addition all one-dimensional $A$-modules are coherent.

The main theorem of \cite{P-coh} asserts that if $E_n$ is ample then the $\Z$-algebra $A$ is coherent and the natural
functor $E\mapsto \oplus_{i<0} \Hom(E_i,E)$ gives an equivalence of categories
\begin{equation}\label{eqcat}
\AA\simeq \cohproj A
\end{equation}
where $\cohproj A$ is the quotient of the category of coherent right $A$-modules by the subcategory of finite-dimensional
modules. We are going to apply this theorem to the category $\CC'$ generated by standard holomorphic bundles on $T=T_{\th,\tau}$.
Recall that in \cite{PS} we identified this category with certain abelian subcategory $\CC^{\th}$ of the
derived category $D^b(X)$ of coherent sheaves on the elliptic curve $X=\C/\Z+\Z\tau$. To define $\CC^{\th}$
one has to consider two full subcategories in the category $\Coh(X)$ of coherent sheaves on $X$: 
$\Coh_{<\th}$ (resp., $\Coh_{>\th}$) is the minimal subcategory of $\Coh(X)$ closed under extensions and containing
all stable bundles of slope $<\th$ (resp., all stable bundles of slope $>\th$ and all torsion sheaves). Then
by the definition
$$\CC^{\th}=\{K\in D^b(X):\ H^{>0}(K)=0, H^0(K)\in\Coh_{>\th}, H^{-1}(K)\in\Coh_{<\th}, H^{<-1}(K)=0\}.$$
Thus, $\CC^{\th}$ contains $\Coh_{>\th}$ and $\Coh_{<\th}[1]$ and these two subcategories generate $\CC^{\th}$
in an appropriate sense. The fact that $\CC^{\th}$ is abelian follows from the torsion theory (see \cite{HRS}). 
Note that the vectors $(\deg(K),\rk(K))\in\Z^2$ for $K\in\CC^{\th}$ are characterized by the inequality
$$\deg(K)-\th\rk(K)>0.$$
In \cite{PS} we showed that a version of Fourier-Mukai transform gives an equivalence $\SS:\CC'\wt{\ra}\CC^{\th}$
(this $\SS$ differs from the transform studied in section 3.3 of \cite{PS} by the shift $K\mapsto K[1]$).
Standard holomorphic bundles correspond under $\SS$ to {\it stable} 
objects of $\CC^{\th}$: the latter are structure sheaves of points and objects
of the form $V[n]$ where $V$ is a stable bundle, $n\in\{0,1\}$.
Moreover, one has
$$\deg\SS(E_{d,c}(\th))=d,\ \rk\SS(E_{d,c}(\th))=-c.$$
It follows that 
$$\rk(\SS(E))=-\deg(E), \ \mu(\SS(E))=\th-\mu(E)^{-1}.$$
The following criterion of ampleness in $\CC'$
is essentially contained in the proof of Theorem 3.5 of \cite{P-real}, where we showed the existence of ample
sequences in $\CC'$.

\begin{thm}\label{amplethm} 
Let $(E_n)$ be a sequence of standard holomorphic bundles on $T$ such that
$\mu(E_n)\to-\infty$ as $n\to-\infty$ and $\rk(E_n)>c$ for all $n<<0$ for some
constant $c>0$. Then $(E_n)$ is an ample sequence in $\CC'$. Moreover, for every
$E\in\CC'$ the natural morphism $\Hom(E_n,E)\otimes_{\C} E_n\ra E$ in $\CC$ is surjective
for $n<<0$.
\end{thm}

\Pf . Let $\FF_n=\SS(E_n)$ be the corresponding sequence of stable objects of $\CC^{\th}$.
Then $\rk(\FF_n)=-\deg(E_n)\to+\infty$ and $\mu(\FF_n)=\th-\mu(E_n)^{-1}\to\th$ as $n\to-\infty$.
Moreover, we have
$$\mu(\FF_n)-\th=\frac{\rk(E_n)}{\rk(\FF_n)}>\frac{c}{\rk(\FF_n)}.$$
Therefore, the same proof as in Theorem 3.5 of \cite{P-real} (where we considered
only the special case $c=1$) shows that
the sequence $(\FF_n)$ is ample in $\CC^{\th}$ and that for every $\FF\in\CC^{\th}$ the morphism
$\Hom(\FF_n,\FF)\otimes_{\C}\FF_n\ra \FF$ is surjective for $n<<0$. Hence, the same assertions hold
for the sequence $(E_n)$ in $\CC'$.
\ed

\begin{thm}\label{cohthm} 
Let $(E_n)$ be a sequence as in Theorem \ref{amplethm} and let
$A=\oplus_{i\le j}\Hom_{\CC}(E_i,E_j)$ be the corresponding $\Z$-algebra. Then for every
holomorphic bundle $E$ on $T$ the $A$-module $M(E)=\oplus_{i<0}\Hom_{\CC}(E_i,E)$ is coherent.
Also, for every sufficiently small $i_0$ the canonical morphism of $A$-modules
$$\Hom_{\CC}(E_{i_0},E)\otimes P_{i_0}\ra M(E)$$
has finite-dimensional cokernel.
\end{thm}

First, we need a criterion for finite generation of modules $M(E)$ with a weaker assumption on $(E_n)$.

\begin{lem}\label{fingenlem}
Let $E$ be a holomorphic bundle on $T$ and let $C=C(E)$ be the corresponding constant from Theorem \ref{vanishthm}.
Let $(E_i)$ be a sequence of standard holomorphic bundles such that $\mu(E_i)\to\mu$ as $i\to-\infty$,
where $\mu\in\R\cup\{-\infty\}$.
Assume that for some $\eps>0$ one has $\rk(E_i)^2(\mu(E_i)-\mu)>1+\eps$ for all $i<<0$ (this condition is
vacuous if $\mu=-\infty$).
Then  

\noindent
(i) for every sufficiently small $i_0\in\Z$ there exists $i_1=i_1(i_0)$ such that 
for all $i<i_1$ there exists a standard holomorphic bundle $F_i$ fitting into the following 
short exact sequence in $\CC'$:
\begin{equation}\label{mainexseq}
0\ra E_i\stackrel{\can}{\ra}\Hom_{\CC}(E_i,E_{i_0})^*\otimes E_{i_0}\ra F_i\ra 0,
\end{equation}
where $\can$ is the canonical morphism.

\noindent
(ii) Assume in addition that $\mu<C$ and if $\mu=-\infty$ then for some $c>0$ one has $\rk(E_i)>c$ for all $i<<0$. 
Then for all sufficiently small $i_0$ there exists $i_1$ such that for $i<i_1$ one has
$\Ext^1(F_i,E)=0$, where $F_i$ is defined by (\ref{mainexseq}). 
Under the same assumptions the $A$-module $M(E)$ is finitely generated.
\end{lem}

\Pf . (i) Let us denote $r_i=\rk(E_i)$, $\mu_i=\mu(E_i)$. 
If $\mu$ is finite then for every sufficiently small $i_0$ one has $r_{i_0}^2(\mu_{i_0}-\mu)>1$. Therefore,
we can find $i_1<i_0$ such that for $i<i_1$ one has
$r_{i_0}^2(\mu_{i_0}-\mu_i)>1$. If $\mu=-\infty$ then we can take any $i_0$ and then still find $i_1<i_0$ such that
for $i<i_1$ the above inequality holds. We are going to construct $F_i$ in this situation. 
Using the equivalence $\SS:\CC'\wt{\ra}\CC^{\th}\sub D^b(X)$
we can first define $\wt{F}_i\in D^b(X)$ from the exact triangle
$$\SS(E_i)\ra\Hom_{\CC}(E_i,E_{i_0})^*\otimes \SS(E_{i_0})\ra \wt{F}_i\ra \SS(E_i)[1].$$
In other words, $\wt{F}_i$ is the image of $E_i$ under the equivalence
$R_{E_{i_0}}:D^b(X)\ra D^b(X)$ given by the 
right twist with respect to $E_{i_0}$ (see \cite{P-real}, sec.2.3, or \cite{ST}; our functor differs from that of \cite{ST}
by a shift).
It follows that $\Hom_{D^b(E)}(\wt{F}_i,\wt{F}_i)\simeq\Hom_{\CC}(E_i,E_i)\simeq\C$, 
so $\wt{F}_i$ is a stable object and either $\wt{F}_i\in\CC^{\th}$ or $\wt{F}_i\in\CC^{\th}[1]$.
To prove that $\wt{F}_i\in\CC^{\th}$ it suffices to check that $\deg(\wt{F}_i)-\th\rk(\wt{F}_i)>0$.
But
$$\deg(\wt{F}_i)-\th\rk(\wt{F}_i)=\chi(E_i,E_{i_0})r_{i_0}-r_i=((\mu_{i_0}-\mu_i)r_{i_0}^2-1)r_i>0$$
by our choice of $i$. Hence, we have $\wt{F}_i\in\CC^{\th}$ and we can set $F_i=\SS^{-1}(\wt{F}_i)$. 

\noindent
(ii) If $\mu$ is finite then we can choose $i_0$ such that $r_{i_0}^2(\mu_{i_0}-\mu)>1+\eps$ and
$\mu_{i_0}+\eps^{-1}(\mu_{i_0}-\mu)<C$. If $\mu=-\infty$ then we choose $i_0$ such that
$\mu_{i_0}+2r_{i_0}^{-2}<C$ (here we use the assumption that $r_i>c$ for $i<<0$). 
In either case for sufficiently small $i$ we have short exact sequence (\ref{mainexseq}).
Applying to it the functor $\Hom(?,E)$ we get the long exact sequence
$$0\ra\Hom_{\CC}(F_i,E)\ra \Hom_{\CC}(E_{i_0},E)\otimes\Hom_{\CC}(E_i,E_{i_0})\ra \Hom_{\CC}(E_i,E)\ra\Ext^1_{\CC}(F_i,E).$$
Thus, vanishing of $\Ext^1_{\CC}(F_i,E)$ for all $i<<0$ would imply that the $A$-module $M(E)$ is finitely generated.
By the definition of the constant $C$ this vanishing would follow from the inequality $\mu(F_i)<C$.
In the case $\mu\neq-\infty$ we have for $i<<0$ 
$$\mu(F_i)=\mu_{i_0}+\frac{\mu_{i_0}-\mu_i}{r_{i_0}^2(\mu_{i_0}-\mu_i)-1}<
\mu_{i_0}+\eps^{-1}(\mu_{i_0}-\mu_i),$$
so the required inequlity follows for $i<<0$ from our choice of $i_0$. In the case $\mu=-\infty$
we can finish the proof similarly using the inequality
$$\mu(F_i)<\mu_{i_0}+\frac{2}{r_{i_0}^2}$$
that holds for all $i<<0$.
\ed

\noindent
{\it Proof of Theorem \ref{cohthm}}.
By Lemma \ref{fingenlem} for any sufficiently small $i_0$ there exists $i_1$ (depending on $i_0$) such that
we have short exact sequence (\ref{mainexseq}) and the induced sequence of $A$-modules
$$0\ra \oplus_{i<i_1} \Hom_{\CC}(F_i,E)\ra \Hom_{\CC}(E_{i_0},E)\otimes (P_{i_0})_{<i_1}\ra M(E)_{<i_1}\ra 0$$
is exact, where for every $A$-module $M=\oplus M_i$ we set $M_{<n}=\oplus_{i<n}M_i$.
This immediately implies the last assertion of the theorem.
Note that the structure of the $A$-module on $\oplus \Hom_{\CC}(F_i,E)$ is defined using the natural isomorphisms
$\Hom_{\CC}(F_i,F_j)\simeq\Hom_{\CC}(E_i,E_j)$ coming from the equality $F_i=R_{E_{i_0}}(E_i)$, where
$R_{E_{i_0}}$ is the right twist with respect to $E_{i_0}$. 
It suffices to prove that the module $M'(E):=\oplus_{i<i_1}\Hom_{\CC}(F_i,E)$ is finitely generated.
Indeed, this would imply that the module $M(E)_{<i_1}$ is finitely presented and hence coherent (since $A$ is coherent), therefore,
the module $M(E)$ is also coherent. To check that $M'(E)$ is finitely generated we will use the criterion of
Lemma \ref{fingenlem} for the sequence $(F_i)$.  
We have 
$$\mu(F_i)\to\mu=\mu_{i_0}+1/r_{i_0}^2$$
as $i\to-\infty$.
Also, 
$$\rk(F_i)^2(\mu(F_i)-\mu)=\frac{((\mu_{i_0}-\mu_i)r_{i_0}^2-1)r_i^2}{r_{i_0}^2}>\frac{((\mu_{i_0}-\mu_i)r_{i_0}^2-1)c^2}{r_{i_0}^2}\to+\infty
$$
as $i\to-\infty$.
Hence, the conditions of Lemma \ref{fingenlem} will be satisfied once we show that $\mu=\mu_{i_0}+1/r_{i_0}^2$ 
can be made smaller than any given constant by an appropriate choice of $i_0$. But this is of course true
since $\mu_{i_0}+1/r_{i_0}^2<\mu_{i_0}+1/c^2$ and $\mu_{i_0}\to-\infty$ as $i_0\to-\infty$.
\ed

\subsection{Proof of Theorem \ref{mainthm}} 

Let us pick a sequence $(E_n)_{n\in\Z}$ of stable holomorphic bundles satisfying conditions of Theorem \ref{amplethm}
(it is easy to see that such a sequence exists, see the proof of Theorem 3.5 in \cite{P-real}).
Let $E$ be a holomorphic bundle on $T$. Then
by Theorem \ref{cohthm} the module $M=\oplus_i\Hom_{\CC}(E_i,E)$ is coherent, hence,
we can consider the object $E'\in\CC'$ corresponding to this module via the equivalence
(\ref{eqcat}). By the definition this means that there is an isomorphism of $A$-modules
\begin{equation}\label{modisom}
M(E')_{<i_0}\simeq M(E)_{<i_0}
\end{equation}
for some $i_0$. By Theorem \ref{cohthm} 
assuming that $i_0$ is small enough we can ensure that for all $i<i_0$ the canonical morphism
$M(E')_i\otimes P_i\ra M(E')$ has finite-dimensional cokernel.
We claim that there exists a morphism $f:E'\ra E$ in $\CC$ that induces the same isomorphism of $A$-modules
$M(E')_{<i_1}\simeq M(E)_{<i_1}$ for some $i_1<i_0$ as the isomorphism (\ref{modisom}).
Indeed, by Theorem \ref{amplethm} we can find a resolution for $E'$ in $\CC'$
of the form
\begin{equation}\label{resolution}
\ldots\ra V_1\otimes E_{n_1}\ra V_0\otimes E_{n_0}\ra E'\ra 0
\end{equation}
where $V_0=\Hom(E_{n_0},E')$ and $n_j<i_0$ for all $j\ge 0$. Using this resolution we can compute $\Hom_{\CC}(E',E)$:
$$\Hom_{\CC}(E',E)\simeq\ker(V_0^*\otimes\Hom_{\CC}(E_{n_0},E)\ra V_1^*\otimes\Hom_{\CC}(E_{n_1},E)).$$
Using isomorphism (\ref{modisom}) we can identify this space with
$$\ker(V_0^*\otimes\Hom_{\CC}(E_{n_0},E')\ra V_1^*\otimes\Hom_{\CC}(E_{n_1},E'))\simeq\Hom_{\CC}(E',E').$$
Thus, we obtain an isomorphism $\Hom_{\CC}(E',E)\simeq\Hom_{\CC}(E',E')$. We define $f\in\Hom_{\CC}(E',E)$ to be the element
corresponding to the identity in $\Hom_{\CC}(E',E')$. Let us check that $f$ induces the same 
isomorphism as (\ref{modisom}) on some truncations of the modules $M(E')$ and $M(E)$. The definition of $f$ implies
that the composition of the induced morphism $f_*:M(E')\ra M(E)$ with the natural morphism
$V_0\otimes P_{n_0}=M(V_0\otimes E_{n_0})\ra M(E')$ coincides with the morphism 
$V_0\otimes P_{n_0}\ra M(E)$ induced by the isomorphism $V_0=\Hom(E_{n_0},E')\simeq\Hom(E_{n_0},E)$ induced by (\ref{modisom}).
Therefore, our claim follows from the fact that the above morphism $V_0\otimes P_{n_0}\ra M(E')$ induces a surjective
morphism on appropriate truncations. 

Thus, we can assume from the beginning that the isomorphism (\ref{modisom}) is induced by a morphism $f:E'\ra E$.
Next, we are going to construct a morphism $g:E'\ra E$
such that $g\circ f=\id_{E'}$. To do this we note that by Serre duality
$\Hom_{\CC}(E,E')\simeq\Ext^1_{\CC}(E',E)^*$ (see Corollary \ref{Sercor}).
Let us make $n_0$ smaller if needed so that $\Ext^1_{\CC}(E_{n_0},E)=\Ext^1_{\CC}(E_{n_0},E')=0$.
Then the space $\Ext^1_{\CC}(E',E)$ can be computed using resolution (\ref{resolution}):
\begin{equation}\label{extisom}
\Ext^1_{\CC}(E',E)\simeq H^1[V_0^*\otimes\Hom_{\CC}(E_{n_0},E)\ra V_1^*\otimes\Hom_{\CC}(E_{n_1},E)\ra 
V_2^*\otimes\Hom_{\CC}(E_{n_2},E)].
\end{equation}
Indeed, let us define $K_1\in\CC'$ from the short exact sequence
$$0\ra K_1\ra V_0\otimes E_{n_0}\ra E'\ra 0,$$
so that we have the following resolution for $K_1$:
$$\ldots\ra V_2\otimes E_{n_2}\ra V_1\otimes E_{n_1}\ra K_1\ra 0.$$
Then the isomorphism (\ref{extisom}) can be derived from the induced exact sequences
$$\Hom_{\CC}(V_0\otimes,E_{n_0},E)\ra\Hom_{\CC}(K_1,E)\ra\Ext^1_{\CC}(E',E)\ra \Ext^1_{\CC}(V_0\otimes E_{n_0},E)=0,$$
$$0\ra\Hom_{\CC}(K_1,E)\ra\Hom_{\CC}(V_1\otimes E_{n_1},E)\ra\Hom_{\CC}(V_2\otimes E_{n_2},E).$$
Using the fact that isomorphism (\ref{extisom}) is functorial in $E$ such that $\Ext^1_{\CC}(E_{n_0},E)=0$ 
we derive that the morphism $\Ext^1_{\CC}(E',E')\ra\Ext^1_{\CC}(E',E)$ induced by $f$ is an isomorphism.
But there is a natural
functional $\phi'\in\Ext^1_{\CC}(E',E')^*$ given by Serre duality. Let $\phi\in\Ext^1_{\CC}(E',E)^*$
be the corresponding functional. The isomorphism $\Ext^1_{\CC}(E',E)^*\wt{\ra}\Hom_{\CC}(E,E')$
maps $\phi$ to some element $g\in\Hom_{\CC}(E,E')$. 
By functoriality of the Serre duality the following diagram is commutative: 
\begin{equation}\label{maindiagram}
\begin{diagram}
\Hom_{\CC}(E,E')&\rTo{\a}&  \Ext^1_{\CC}(E',E)^*\nonumber\\
\dTo{f^*} &&\dTo^{f^*}\\
\Hom_{\CC}(E',E')&\rTo{\a'}& \Ext^1_{\CC}(E',E')^*
\end{diagram}
\end{equation}
where the vertical arrows are induced by $f$. 
Since $\phi'=\a'(\id_{E'})$, $f^*(\phi)=\phi'$ and $\a(g)=\phi$ we deduce
that $f^*(g)=\id_{E'}$, i.e. $g\circ f=\id_{E'}$.
Therefore, we have $E\simeq E'\oplus E''$ for some holomorphic bundle $E''$
such that $\Hom_{\CC}(E_i,E'')=0$ for $i<i_0$. But Theorem \ref{vanishthm} implies that
$\Ext^1_{\CC}(E_i,E'')=0$ for all sufficiently negative $i$. Together
with Corollary \ref{RRcor} this implies that $E''=0$.
\ed

\end{document}